\documentclass[leqno,12pt]{amsart}
\usepackage{amsmath} 
\usepackage{amssymb,amsxtra,amsthm,amscd}
\usepackage{graphicx} 
\usepackage[pdftex]{hyperref}
\usepackage[margin=34mm]{geometry}

\usepackage{blaom} 
\DeclareMathOperator{\rad}{rad}
\usepackage{pdfsync}
\begin{document}
\title%
[Lie algebroids]%
{The Lie algebroid associated with a hypersurface} \author{Anthony
  D.~Blaom}
\address{E-mail: {\tt anthony.blaom@gmail.com}}%
\thispagestyle{empty}
\begin{abstract} In this note we motivate the definition and use of
  Lie algebroids by revisiting the problem of reconstructing a
  hypersurface in Euclidean space from infinitesimal data.
\end{abstract}
\maketitle


\section{Introduction}\label{sec1}
Lie algebroids are the infinitesimal version of Ehresmann's {\it
  categories differentiables} \cite{Ehresmann_59}, now called Lie
groupoids. They were first introduced and studied by J. Pradines in
1967 \cite{Pradines_67}. Twenty years later Alan Weinstein and others
pointed out important connections with Poisson geometry
\cite{Weinstein_87,Coste_etal_87}. There was a resurgence of interest,
deep results followed, and the subject remains very active.
Applications include foliation theory \cite{Moerdijk_Mrcun_03},
noncommutative geometry \cite{CannasdaSilva_Weinstein_99,Landsman_06},
Lagrangian mechanics and control theory \cite{Cortes_etal_06},
elasticity theory \cite{Jimenez_etal_2016}, and complex analysis
\cite{Gualtieri_etal_13}.

Lie algebroids also arise naturally in elementary differential
geometry. The purpose of the present note is to give a simple example
of this by relating the recent observations of \cite{Blaom_F,Blaom_G}
in a special case. Our note is meant as an invitation to Lie
algebroids, about which no prior knowledge is assumed. Proper
introductions and further references can be found in
\cite{CannasdaSilva_Weinstein_99,Crainic_Fernandes_11,Dufour_Nguyen_05,Mackenzie_03}. For
detailed historical notes, prior to 2003, see \cite{Mackenzie_03}. We
assume familiarity with basic calculus on manifolds and Lie groups.

\section{On the Killing fields tangent to a hypersurface}\label{sec2}
The infinitesimal isometries of ${\mathbb R}^{n+1} $ are the {\df
  Killing fields}, those vector fields whose local flows are rigid
motions. To characterise a hypersurface $\Sigma $ in
${\mathbb R}^{n+1}$ up to rigid motions, it suffices to understand the
relationship between $\Sigma $ and the Lie algebra ${\mathfrak g} $ of
Killing fields.\footnote{For simplicity we shall behave as if
  $\Sigma $ is embedded in ${\mathbb R}^{n+1} $. However, with obvious
  modifications, our constructions apply to hypersurfaces that are
  merely immersed.} To this end, consider those smooth functions
$\mathbf{X} \colon \Sigma \rightarrow {\mathfrak g} $ on $\Sigma $ that are
{\df tangent} to $\Sigma$, in the sense that
\begin{equation*}
  \mathbf{X}(x)(x) \in T_x \Sigma,
\end{equation*}
for all $x \in \Sigma $. Such functions amount to sections of the
vector bundle ${A}$ over $\Sigma $ defined by
\begin{equation*}
  {A} = \{(X,x) \in {\mathfrak g} \times \Sigma \suchthat X(x) \in
  T_x \Sigma\}.
\end{equation*}
Evidently, $A$ is a subbundle of the trivial bundle
${\mathfrak g} \times \Sigma $ of rank $n(n+3)/2$.

By construction, we have a canonical vector bundle morphism
$\# \colon {A} \rightarrow T \Sigma $ defined by
$\#(X,x)=X(x)$. This map happens to be surjective, suggesting that
we think of ${A}$ as a `thickening' of $T\Sigma $. Also, if
${\mathfrak h} \subset A$ denotes the kernel of $\#$ then the fibre
${\mathfrak h}|_x$ can be identified with the subalgebra of Killing
fields vanishing at $x \in \Sigma $. In particular, ${\mathfrak h} $
is a bundle of Lie algebras.

Now $A$ itself is {\em not} a bundle of Lie algebras. However, we have
the following intriguing observation: Just as the Jacobi-Lie bracket
on vector fields on $\Sigma $ makes the space of vector fields
$\Gamma(T \Sigma)$ into a Lie algebra, and satisfies the Leibniz
identity,
\begin{equation*}
  [X,fY]=f[X,Y]+df(X)Y, \qquad X,Y \in \Gamma(T\Sigma),
\end{equation*}
for all smooth functions $f \colon \Sigma \rightarrow {\mathbb R} $,
so the section space $\Gamma ({A})$ has a god-given Lie algebra
bracket satisfying an almost identical identity, namely
\begin{equation}
  [X,fY]=f[X,Y]+df(\#X)Y, \qquad X,Y \in \Gamma({A}).\label{wert}
\end{equation}
Note here that we need the map $\# \colon {A} \rightarrow T\Sigma $
for we do not otherwise have a way of differentiating the function $f$
along a section of ${A}$.

How is this bracket on ${A}$ defined? First, define a bracket on
${\mathfrak g} \times \Sigma $ in the most naive way:
\begin{equation*}
  \{X,Y\}(x)=[X(x),Y(x)]_{\mathfrak g}.
\end{equation*}
Here we are viewing sections of
$A \subset {\mathfrak g} \times \Sigma $ as ${\mathfrak g} $-valued
functions, and the bracket on the right is the bracket on
${\mathfrak g} $, the restriction of the Jacobi-Lie braket to the
subalgebra of Killing fields. This bracket does not satisfy a
Leibniz-type identity; rather it is `algebraic' (bilinear over the
ring of all smooth functions on $\Sigma $). However, we have:
\begin{proposition}
  Let $\nabla $ denote the canonical flat connection on the trivial
  bundle ${\mathfrak g} \times \Sigma $. Then the bracket on
  $\Gamma({A})$ given by 
  \begin{equation*}
    [X,Y]=\nabla_{\#X}Y-\nabla_{\#Y}X +\{X,Y\}
  \end{equation*}
  is well-defined and satisfies the Leibniz identity
  \eqref{wert}. Moreover, with respect to this bracket, the map
  $\Gamma({A}) \rightarrow \Gamma(T \Sigma )$ induced by
  $\# \colon {A} \rightarrow T \Sigma $ is a Lie algebra homomorphism.
\end{proposition}
It is not immediately obvious that the bracket is well-defined as
${A} \subset {\mathfrak g} \times \Sigma $ is not $\nabla $-invariant,
unless $\Sigma $ has lots of symmetry (e.g., is a hyperplane or
hypersphere). Moreover, $\Gamma(A)$ is not generally closed under the
algebraic bracket $\{\,\cdot\,,\,\cdot\,\}$. The reader will find it
instructive to prove this proposition on her own. A proof is included
in the appendix.

\section{Lie algebroids}
\begin{definition}
  A {\df Lie algebroid} with {\df base} $\Sigma $ is a vector bundle
  $A$ over $\Sigma $, together with a vector bundle morphism
  $\# \colon A \rightarrow T \Sigma $ covering the identity and called
  the {\df anchor}, and a Lie bracket on its space of sections, such
  that the Leibniz identity \eqref{wert} holds.
\end{definition}
Every tangent bundle is a Lie algebroid with the identity map as
anchor, and every ordinary Lie algebra is a Lie algebroid over a
point. Every Lie groupoid (see Section \ref{orange}) has a Lie
algebroid as its infinitesimalization. Other important examples are
recalled in the appendix.

To explain the significance of the Lie algebroid ${A}$ constructed
above for a hypersurface $\Sigma \subset {\mathbb R}^{n+1} $ we must
say something about {\df morphisms} in the category of Lie algebroids.
If a vector bundle map $\omega \colon A_1 \rightarrow A_2$ between Lie
algebroids covers the identity on a common base $\Sigma $, then we
declare $\omega $ to be a Lie algebroid morphism if the corresponding
map of section spaces $\Gamma(A_1) \rightarrow \Gamma(A_2)$ is a Lie
algebra homomorphism. In particular, the anchor
$\# \colon A \rightarrow T \Sigma $ of a Lie algebroid is always a Lie
algebroid morphism, a corollary of the Jacobi and Leibniz identities;
for a proof see the appendix. When the base of $A_1$ and $A_2$ are
different then the definition of a morphism is more complicated and we
give it here in a special case only: Suppose $A$ is a Lie algebroid
over $\Sigma $ and ${\mathfrak g} $ a Lie algebroid over a {\em point}
(a Lie algebra). Then a linear map
$\omega \colon A \rightarrow {\mathfrak g} $ is a Lie algebroid
morphism if the obvious extension to a morphism of vector bundles
$\bar\omega \colon A \rightarrow {\mathfrak g} \times \Sigma $ satisfies
\begin{equation}
  \bar\omega ([X,Y])=\nabla_{\#X} \bar\omega (Y)-\nabla_{\#Y}\bar\omega (X)
  +\{\bar\omega (X),\bar\omega (Y)\}.\mathlab{huy}
\end{equation}
Here $\nabla $ and $\{\,\cdot\,,\,\cdot\,\}$ have the meanings given
in Section \ref{sec1}. (In particular, a Lie algebroid morphism
$\omega \colon T \Sigma \rightarrow {\mathfrak g} $ is the same thing
as a Mauer-Cartan form on $\Sigma $, and \eqref{huy} generalizes the
classical Mauer-Cartan equations \cite{Sharpe_97}.)

\section{An infinitesimal characterization of
  hypersurfaces} \label{orange} \setcounter{equation}{0}%
According to the preceding definition, if $A$ is the Lie algebroid
associated with a hypersurface $\Sigma \subset {\mathbb R}^{n+1} $,
then the composite $\omega \colon A \rightarrow {\mathfrak g} $ of the
the inclusion $A \subset {\mathfrak g} \times \Sigma $ with the
projection ${\mathfrak g} \times \Sigma \rightarrow {\mathfrak g} $ is
a Lie algebroid morphism into the Lie algebra ${\mathfrak g} $ of
Killing fields. We call the map
$\omega \colon A \rightarrow {\mathfrak g} $ the {\df logarithmic
  derivative} of the embedding $\Sigma \subset {\mathbb R}^{n+1} $,
for it may be viewed as a generalization of \'Elie Cartan's
logarithmic derivative of a smooth map into a Lie group (see, e.g.,
\cite{Sharpe_97}) known also as its Darboux derivative. According to
the following theorem, any hypersurface can be reconstructed from its
logarithmic derivative. In the statement
$\rad {\mathfrak g} \subset {\mathfrak g} $ denotes the subalgebra of
constant vector fields (the radical of the Lie algebra of Killing
fields).
\begin{theorem}[The abstract Bonnet theorem]
  Let $A $ be any Lie algebroid over a smooth, orientable,
  simply-connected, $n$-dimensional manifold $\Sigma$, and let
  $\omega \colon A \rightarrow {\mathfrak g} $ be a Lie algebroid
  morphism, where ${\mathfrak g} $ is the Lie algebra of Killing
  fields on ${\mathbb R}^{n+1} $. Assume:
  \begin{conditions}
  \item $\rank A = n(n+3)/2$. \label{one}
  \item The anchor $\# \colon A \rightarrow T \Sigma $ is surjective
    ({\df transitivity}).\label{two}
  \item $\omega \colon A \rightarrow {\mathfrak g} $ is injective on
    fibres.\label{twob}
  \item For some (and consequently any) $x_0 \in \Sigma $, $\omega $
    maps the intersection of the kernel of $\#$ with the fibre
    $A|_{x_0}$ to a subspace of $\mathfrak g$ transverse to
    $\rad {\mathfrak g}$ ({\df transversality}).\label{three}
  \end{conditions}
  Then $\Sigma $ can be realised as an immersed hypersurface in
  ${\mathbb R}^{n+1} $ in such a way that its associated Lie algebroid is
  isomorphic to $A $, and such that its logarithmic derivative is
  $ \omega \colon A \rightarrow {\mathfrak g}$.
\end{theorem}

\begin{remark}
  If a submanifold $\Sigma \subset {\mathbb R}^{n+1}$ is replaced by
  its image $\phi(\Sigma)$ under some rigid motion
  $\phi \colon {\mathbb R}^{n+1} \rightarrow {\mathbb R}^{n+1} $, then
  it logarithmic derivative is altered, but in a predictable way: just
  compose the old logarithmic derivative with
  $\Adjoint_\phi \colon {\mathfrak g} \rightarrow {\mathfrak g} $,
  where $\Adjoint_\phi \xi $ is the pushforward of the vector field
  $\xi $ under $\phi $. In this way, we may regard the logarithmic
  derivative of an embedding as an invariant under rigid
  motions. Moreover, up to rigid motions, every embedding can be
  reconstructed from this invariant, by the preceding theorem.
\end{remark}

Noting that the sufficient conditions \eqref{one}--\eqref{three} are
also necessary, we now sketch a proof of the abstract Bonnet theorem
under the additional global assumption that $A$ is {\df integrable},
by which we mean that $A$ is the infinitesimalization of some Lie
groupoid ${\mathcal G} $.  This extra assumption is also necessary but
turns out to be redundant \cite{Blaom_F}.  In the proof we describe
Lie groupoids sufficiently that the magnanimous reader will grasp the
main ideas, and appreciate how Lie theory can be applied to problems
of differential geometry in novel ways. In particular, let us
emphasise that the same proof more-or-less delivers a characterization
of arbitrary smooth maps $f \colon \Sigma \rightarrow G/H$ of a
simply-connected manifold $\Sigma $ into an arbitrary homogeneous
space $G/H$ (`Klein geometry'); if $\Sigma $ is not simply-connected,
then there is a global obstruction to reconstructing such maps from
infinitesimal data called the {\df monodromy}, also naturally
formulated using Lie algebroid language \cite{Blaom_F}. Familiarity
with our proof is not required in the remainder of this note.

\begin{proof}[Sketch of proof of the theorem]
  A Lie groupoid is a Lie group ${\mathcal G} $ whose identity element
  has been `blown up' into a smooth embedded submanifold $\Sigma $
  (the {\df base}) and for which multiplication is consequently
  defined only partially. In more detail, we should understand that
  ${\mathcal G} $ is a smooth manifold whose points represent the
  arrows (morphisms) of a category, each arrow $g \in {\mathcal G} $
  beginning at one identity element $\alpha(g) \in \Sigma $ and ending
  at another, $\beta(g) \in \Sigma $.  That is, the identity elements
  represent the objects of the category, and at the same time
  represent the identity arrows. Two arrows
  $g_1,g_2 \in {\mathcal G} $ are {\df multipliable} if they can be
  composed as arrows of the category, and their product is, by
  definition, their composition in the category. The partial
  multiplication is assumed to be smooth, as are the {\df source} and
  {\df target maps}
  $\alpha, \beta \colon {\mathcal G} \rightarrow \Sigma $ of the
  groupoid, which are additionally required to be submersions. Every
  element of ${\mathcal G} $ is assumed to have a multiplicative
  inverse. A morphism $\Omega $ between Lie groupoids is a smooth
  functor of the underlying categories. In particular, this means
  \begin{equation}
    \Omega (gh)=\Omega (g) \Omega (h),\label{yellow}
  \end{equation}
  whenever the product $gh $ is defined.

  For the geometer, the canonical example of a Lie groupoid is the
  {\df orthonormal frame groupoid} of a Riemannian manifold $\Sigma$,
  whose arrows consist of the isometries between possibly different
  tangent spaces of $\Sigma $, the base of the groupoid.

  The infinitesimalization of a Lie groupoid is a Lie algebroid. Just
  as the Lie algebra of a Lie group is the tangent space to the
  identity, so more generally, the Lie algebroid of a Lie groupoid
  ${\mathcal G} $ is the normal vector bundle of the base $\Sigma $ of
  identity elements in ${\mathcal G} $. We omit here the description
  of the bracket and anchor.

  Now suppose the Lie algebroid $A$ in the theorem is the
  infinitesimalization of a Lie groupoid ${\mathcal G} $ over
  $\Sigma $. In analogy with the case of Lie groups, we may suppose
  that ${\mathcal G} $ is `simply-connected' which here means the
  source projection $\alpha \colon {\mathcal G} \rightarrow \Sigma $,
  taking an arrow to its starting point, has simply-connected
  fibres. The theorem known as Lie II states that every morphism of
  Lie algebras ${\mathfrak g}_1 \rightarrow {\mathfrak g}_2$ is the
  infinitesimalization of a homomorphism $G_1 \rightarrow G_2$ between
  Lie groups having ${\mathfrak g}_1$ and ${\mathfrak g}_2$ as their
  infinitesimalizations, assuming $G_1$ is simply-connected.  Lie II
  generalizes to Lie groupoids, so that the morphism
  $\omega \colon A \rightarrow {\mathfrak g} $ lifts to Lie groupoid
  morphism $\Omega \colon {\mathcal G} \rightarrow G$, where $G $
  is the isometry group of ${\mathbb R}^{n+1}$.

  To construct an immersion
  $f \colon \Sigma \rightarrow {\mathbb R}^{n+1} $ using $\Omega$ we
  first choose an appropriate target $m_0 \in {\mathbb R}^{n+1} $ for
  the point $x_0 \in \Sigma $ appearing in the transversality
  condition \eqref{three}. To do this, let
  ${\mathfrak h}_{x_0} \subset A$ denote the intersection of the
  kernel of $\# \colon A \rightarrow T\Sigma $ with $A|_{x_0}$. Then,
  by transversality and a dimension count,
  $\omega ({\mathfrak h}_{x_0})$ is a complement of
  $\rad {\mathfrak g} $ in ${\mathfrak g} $. Every such complement is
  the subalgebra ${\mathfrak g}_{m_0}\subset {\mathfrak g} $ of all
  Killing fields vanishing at some point $m_0 \in {\mathbb R}^{n+1} $,
  so that $ \omega({\mathfrak h}_{x_0})= {\mathfrak g}_{m_0}$ and, in
  particular,
  \begin{equation}
    \omega({\mathfrak h}_{x_0}) \subset {\mathfrak g}_{m_0}.\label{grab}
  \end{equation}
  We define 
  \begin{equation*}
    f(x)=\Omega(g) \cdot m_0,
  \end{equation*}
  where $g \in {\mathcal G} $ is any arrow from $x_0$ to $x$. Such an
  arrow exists, by what is called the {\df transitivity} of
  ${\mathcal G} $, following from the transitivity of $A $. The main
  issue is to show that $f(x)$ is independent of the choice of arrow
  $g$. From \eqref{yellow} it suffices to show that
  \begin{equation}
    \Omega(H_{x_0})\subset G_{m_0},\label{grabs}
  \end{equation}
  where $H_{x_0}\subset {\mathcal G} $ is the group of all arrows
  beginning and ending at $x_0$, and $G_{m_0}$ is the group of
  isometries fixing $m_0$. It will not surprise the reader to learn
  that $H_{x_0}$ is a Lie group whose infinitesimalization is the Lie
  algebra ${\mathfrak h}_{x_0}$, so that condition \eqref{grab} is
  precisely the infinitesimalization of condition \eqref{grabs}.  To
  show that the former implies the latter, it accordingly suffices to
  show that $H_{x_0}$ is connected. However, it turns out that the
  restriction of the target projection
  $\beta \colon {\mathcal G} \rightarrow \Sigma $ to the connected
  fibre $\alpha^{-1}(x_0)$, is a principal $H_{x_0}$-bundle over
  $\Sigma $.  The simple-connectivity of $\Sigma $ and the long exact
  sequence in homotopy for the principal bundle shows that $H_{x_0}$
  is indeed connected.

  To prove that $f \colon \Sigma \rightarrow {\mathbb R}^n $ is an
  immersion with the desired properties requires that we describe the
  infinitesimalization functor from Lie groupoids to Lie algebroids in
  more detail and this not attempted here.
\end{proof}

\section{The first and second fundamental forms}
We suppose the reader is already acquainted with the classical Bonnet
theorem for hypersurfaces \cite{Tenenblat_71,Jacobowitz_82}, also
known as the fundamental theorem for hypersurfaces, and is no-doubt
wondering: Where is the second fundamental form? Where are the
Gauss-Codazzi equations?  Naturally, this information is encoded in
the logarithmic derivative of the hypersurface.

  Let $\omega \colon A \rightarrow {\mathfrak g} $ be any Lie algebroid
  morphism satisfying the hypotheses of the abstract Bonnet
  theorem above. This data amounts to extra structure on the two-manifold
  $\Sigma $. We now define two symmetric tensors $g_\omega$ and
  $\ii_\omega $ which coincide with the first and second fundamental
  forms of the hypersurface when
  $\omega \colon A \rightarrow {\mathfrak g} $ is the logarithmic
  derivative of some embedding
  $\Sigma \hookrightarrow {\mathbb R}^{n+1} $.

  On account of condition \eqref{twob} of the theorem, we may regard
  $A$ as a subbundle of $E := {\mathfrak g} \times \Sigma $ and do so
  from now on. Let ${\mathfrak h} \subset A$ denote the kernel of the
  anchor $\# \colon A \rightarrow T \Sigma $. Notice that when $A $ is
  the Lie algebroid associated with an embedded hypersurface
  $\Sigma $, we may evidently identify the quotient $E/{\mathfrak h} $
  with $T_\Sigma {\mathbb R}^{n+1}$, the pullback of the tangent
  bundle of ${\mathbb R}^{n+1} $ to $\Sigma $, which contains
  $T \Sigma $ as a subbundle. This follows from the fact that every
  tangent space of ${\mathbb R}^{n+1}$ is spanned by the Killing
  fields evaluated there. In the general case, we may continue to
  regard $T \Sigma $ as a subbundle of $E/{\mathfrak h} $, for the
  transitivity of $A$, condition \eqref{two}, means we may regard the
  injection $A/{\mathfrak h} \hookrightarrow E/{\mathfrak h} $ as a
  map $T \Sigma \hookrightarrow E/{\mathfrak h} $. However, conditions
  \eqref{one} and \eqref{three} of the theorem ensure that the
  projection $E \rightarrow E/{\mathfrak h} $ maps
  $\rad {\mathfrak g} \times \Sigma $ isomorphically onto
  $E/{\mathfrak h} $. So we in fact have an inclusion
  $T\Sigma \hookrightarrow \rad {\mathfrak g} \times \Sigma $. Making
  the obvious identification
  $\rad {\mathfrak g} \cong {\mathbb R}^{n+1} $, we recover an
  inclusion $\iota \colon T_x\Sigma \hookrightarrow {\mathbb R}^{n+1}$
  (for each $x \in \Sigma $) as in the case of a bona fide
  immersion. We now mimic the usual construction of the first and
  second fundamental forms for immersions, defining\footnote{A
    slightly different but equivalent definition is given in
    \cite{Blaom_G}.}
  \begin{align*}
    g_\omega(X,Y)&=\langle\!\langle \iota X, \iota Y\rangle\!\rangle,\\
    \ii_\omega(X,Y)&=-\langle\!\langle \iota X,d{\mathbf n}(Y)\rangle\!\rangle,
  \end{align*}
  where the $\langle\!\langle \,\cdot\,,\,\cdot\,\rangle\!\rangle$ is
  the standard inner product on ${\mathbb R}^{n+1} $. Here
  ${\mathbf n} \colon \Sigma \rightarrow S^n \subset {\mathbb R}^{n+1}
  $
  is the analogue of the usual Gauss map, defined by declaring
  ${\mathbf n}(x)$ to be orthogonal to
  $\iota(T_x \Sigma) \subset{\mathbb R}^{n+1} $ and have length
  one. The ambiguity in its sign is resolved by supposing $\Sigma $
  has been oriented.  The reader will now readily establish the
  following:
\begin{addendum}
  If $\Sigma $ is oriented and
  $\omega \colon A \rightarrow {\mathfrak g} $ is the logarithmic
  derivative of an immersion $\Sigma \subset {\mathbb R}^{n+1} $, then
  $g_\omega $ and $\ii_\omega $ coincide with the usual first and
  second fundamental forms. In particular, the first and second
  fundamental forms of the realisation in the abstract Bonnet theorem
  will be $g_\omega $ and $\ii_\omega $.
\end{addendum}

The classical Bonnet theorem follows from the fact that for any metric
$g $ and symmetric tensor $\ii$ on $\Sigma $ satisfying the
Gauss-Codazzi equations, we can construct a Lie algebroid morphism
$\omega \colon A \rightarrow {\mathfrak g} $ satisfying the hypotheses
of the abstract Bonnet theorem above, and such that $g_\omega=g$ and
$\ii_\omega =\ii$. For details, see \cite{Blaom_G}.

\appendix
\section{Appendix}
\subsection*{Examples of Lie algebroids} 
We now mention some of the more important examples of Lie
algebroids. If a Lie algebra ${\mathfrak g} $ acts on a manifold
$\Sigma $, then the trivial bundle ${\mathfrak g} \times \Sigma $
becomes a Lie algebroid over $\Sigma $ called an {\df action
  algebroid}, whose anchor is the action map (see below). The image of
the anchor map is tangent to the orbits of the action in this case;
more generally an anchor has as image an involutive $n$-plane field
(possibly with singularities) whose integrating leaves are called the
{\df orbits} of the Lie algebroid. The $n$-plane field tangent to a
regular foliation on a manifold $\Sigma$ is a subalgebroid of the
tangent bundle of $\Sigma$ whose orbits are the leaves of the
foliation. If $\Sigma $ is a Poisson manifold, then $T^* \Sigma $ is a
Lie algebroid whose orbits are the symplectic leaves and, conversely,
for every Lie algebroid $A$, the dual bundle $A^*$ is a (linear)
Poisson manifold. 

If $P$ is a principal $H$-bundle over $\Sigma $, then the quotient
$TP/H$ is a transitive Lie algebroid over $\Sigma $ known as the {\df
  Atiyah Lie algebroid}. In particular, every $G$-structure determines
a Lie algebroid, and every classical Cartan geometry determines a
Koszul connection on the associated Atiyah Lie algebroid.

If $L$ is the canonical line bundle of a contact structure, then its
first jet bundle $J^1L$ is a Lie algebroid. If $A$ is a Lie algebroid,
then so is every jet bundle $J^k A$. 

For more details, more examples, and generalizations of the examples
just given, we refer the reader to the introductory treatments already
cited, the references contained therein, and the following incomplete
list of works:
\cite{Courant_90,Hitchin_03,Blaom_05,Blaom_16b,Cortes_etal_06,Blaom_12,Gualtieri_etal_13,Crainic_Salazar_15}.

The global counterpart of a Lie algebroid is called a {\df Lie
  groupoid}. Much of the relationship between Lie groups and Lie
algebras carries over to the `oid' case, with one important exception:
Not every Lie algebroid is the infinitesimalization of a Lie groupoid,
and the obstructions to `integrating' a Lie algebroid to a Lie
groupoid are subtle \cite{Crainic_Fernandes_03}.

\subsection*{Anchor maps are morphisms}
Here is a proof of the fact that the anchor $\# \colon A \rightarrow
T\Sigma $ of a Lie algebroid is a Lie algebroid morphism. Let $X,Y,Z$
be arbitrary sections of $A $ and $f \colon \Sigma \rightarrow
{\mathbb R} $ a smooth function. Then the Jacobi identity for the
Lie bracket  on $\Gamma(A)$ implies 
\begin{equation}
  [[X,Y],fZ]+[[fZ,X], Y] + [[Y,fZ],X]=0.\label{bone}
\end{equation}
On the other hand, the Leibniz identity implies
\begin{align*}
  [[X,Y],fZ]&=f[[X,Y],Z] + {\mathcal L}_{\#[X,Y]}f\,Z,\\
  [[fZ,X], Y]&=f[Y,[X,Z]] + {\mathcal L}_{\#Y}f\,[X,Z] + {\mathcal
               L}_{\#Y}{\mathcal L}_{\#X}f\,Z,\\
  [[Y,fZ],X]&=-f[X,[Y,Z]]-{\mathcal L}_{\#X}f\,[Y,Z]-{\mathcal
              L}_{\#Y}f\,[X,Z] -{\mathcal L}_{\#X}{\mathcal L}_{\#Y}f\,Z,
\end{align*}
where ${\mathcal L} $ denotes Lie derivative.  Substituting these
equations into \eqref{bone} delivers
\begin{align*}
  {\mathcal L}_{\#[X,Y]}f\, Z + {\mathcal L}_{\#Y}{\mathcal
  L}_{\#X}f\,Z - {\mathcal L}_{\#X}{\mathcal L}_{\#Y}f\,Z&=0\\
  \text{i.e.,}\enspace {\mathcal L}_{\#[X,Y]-[\#X,\#Y]}f\, Z&=0,
\end{align*}
the second line following from the definition of the Jacobi-Lie
bracket. Since $f$ and $Z$ are arbitrary, we conclude
\begin{equation*}
  \#[X,Y]=[\#X,\#Y]; \qquad X,Y \in \Gamma(A). 
\end{equation*}

\subsection*{Action algebroids}
Now suppose a Lie algebra ${\mathfrak g} $ acts on a manifold $M$, and
let
$\xi \mapsto \xi^\dagger \colon {\mathfrak g} \rightarrow \Gamma(TM)$
denote the corresponding Lie algebra homomorphism. For example, we may
take $M={\mathbb R}^{n+1}$ and let ${\mathfrak g} $ be the Killing
fields on ${\mathbb R}^{n+1} $. A bracket on sections of the trivial
bundle ${\mathfrak g} \times M$ is defined by
\begin{equation*}
      [X,Y]=\nabla_{\#X}Y-\nabla_{\#Y}X +\{X,Y\},
\end{equation*}
where $\nabla $ is the canonical flat connection on
${\mathfrak g} \times M$ and $\{X,Y\}(m):=[X(m),Y(m)]_{\mathfrak g} $,
$m \in M$. Since the algebraic bracket $\{\,\cdot\,,\,\cdot\,\}$ is
$\nabla$-parallel, it is easy to see that $[\,\cdot\,,\,\cdot\,]$
satisfies the Jacobi identity. The vector bundle
${\mathfrak g} \times M$ becomes the {\df action algebroid} associated
with the Lie algebra action if we define the anchor
$\# \colon {\mathfrak g} \times M \rightarrow TM$ to be the action map
$\#(\xi,m)=\xi^\dagger(m)$.

In Lie algebroid jargon, the Lie algebroid $A$ associated with a
hypersurface $\Sigma \subset {\mathbb R}^{n+1} $ is the pullback, in
the category of Lie algebroids, under the immersion
$\Sigma \hookrightarrow {\mathbb R}^{n+1}$, of the action algebroid
${\mathfrak g} \times {\mathbb R}^{n+1}$. Pullbacks are defined in
\cite[\S4.2]{Mackenzie_05}.

\subsection*{Proof of Proposition \ref{sec2}}
Consider the action of ${\mathfrak g} $ on ${\mathbb R}^{n+1}$. A
simple proof of the proposition rests on the fact that the action map
$\# \colon {\mathfrak g} \times {\mathbb R}^{n+1} \rightarrow T
{\mathbb R}^{n+1} $
is a Lie algebroid morphism, since it is the anchor map of the action
algebroid ${\mathfrak g} \times {\mathbb R}^{n+1} $. Whence,
  \begin{equation}
  \#[\tilde X,\tilde Y]=[\#\tilde X,\#\tilde Y];\qquad \tilde X,\tilde
  Y \in \Gamma({\mathfrak g} \times  {\mathbb R}^{n+1} ).\mathlab{rabbit}
\end{equation}
So let $X,Y$ be sections of $A \subset {\mathfrak g} \times \Sigma $.
Extend $X$ and $Y$ to sections $\tilde X$ and $\tilde Y$ of
${\mathfrak g} \times {\mathbb R}^{n+1} $. Since $\tilde X$ and
$\tilde Y$ restricted to $\Sigma$ are sections of $A$, the vector
fields $\#\tilde X$ and $\#\tilde Y$ on ${\mathbb R}^{n+1} $ are
tangent to $\Sigma $. In particular, for any $x \in \Sigma $, we have
\begin{equation}
  [\#\tilde X, \# \tilde Y](x) \in T_x \Sigma.\mathlab{got}
\end{equation}
But, on account of \eqref{rabbit}, the left-hand side of \eqref{got}
coincides with $\#[\tilde X, \tilde Y](x)=\#[X,Y](x)$, which shows
$[X,Y]$ is a section of $A$.

\end{document}